\documentclass{siamltex}
\usepackage{float}
\usepackage{geometry}                
\usepackage{graphicx}
\usepackage{amssymb}
\usepackage{amsmath}
\usepackage{algorithmic}
\usepackage{algorithm}
\usepackage{subfigure}
\newtheorem{remark}[theorem]{Remark}
\newcommand{\ha}{\hat{e}}
\newcommand{\hr}{\hat{r}}
\newcommand{\be}{\begin{equation}}
\newcommand{\ee}{\end{equation}}
\newcommand{\ep}{e^{ikx}}
\newcommand{\en}{e^{-ikx}}
\newcommand{\hL}{\hat{L}}
\newcommand{\hap}{\ha_{+}}
\newcommand{\han}{\ha_{-}}
\newcommand{\hanp}{\ha_{\pm}}
\newcommand{\hrp}{\hr_{+}}
\newcommand{\hrn}{\hr_{-}}
\newcommand{\hLp}{\hL_{+}}
\newcommand{\hLn}{\hL_{-}}

\newcommand{\bx}{\bar{x}}
\title{Wave-ray  algorithms for  Helmholtz equations with variable wave numbers: a one-dimensional implementation of  two-dimensional ideas}
\author{I.~Livshits\thanks{Department of Mathematical Sciences, Ball State University, Muncie, IN 47304, e-mail {\sl ilivshits@bsu.edu}}}
\begin{document}
\maketitle

\begin{abstract}
The subject of this paper is multigrid solvers for  Helmholtz operators with large wave numbers. Algorithms presented here are  variations of the {\sl wave-ray} solver which is modified to  allow efficient solutions for operators with constant,  continuous, and discontinuous wave numbers.  
Both geometric and algebraic multigrid frameworks are employed, all yielding efficient and scalable solvers at very little additional costs, compared to standard multigrid techniques.  The algorithms are implemented in one dimension, but with a clear extension to higher dimensions.

\end{abstract}
\section{Introduction}

We consider  a one-dimensional Helmholtz equation 
\be
Lu(x) = \Delta u(x) + k^2(x) u(x) = f(x), \quad x \in \Omega, 
\label{eq:helmholtz}
\ee
accompanied  by the first order Sommerfeld boundary conditions. For  $\Omega = [a, b]$ they read  as 
\be 
\begin{aligned}
& \biggl(\frac{du(x)}{dx}+ik(x) u(x)\biggr)_{|_{x = a}} & = 0, \\
& \biggl(\frac{du(x)}{dx}-ik(x) u(x)\biggr)_{|_{x = b}} & = 0. 
\end{aligned}
\label{eq:bc}
\ee

When discretized on a sufficiently fine grid with mesh-size $h$ (see e.g. \cite{BL06} for what {\sl sufficient}  means) the differential problem is reduced to a system of linear equations 
\be 
A^h u^h = f^h,  \quad  A^h \in \mathbf{C}^{n\times n}, \quad f^h\in \mathbf{C}^n.
\label{eq:discrete}
\ee 
Finding $u^h \in \mathbf{C}^n$  is the  goal of our numerical solvers.

\begin{remark}
 Solving a  one-dimensional Helmholtz equation is  generally  a  significant simplification compared to the application-driven higher dimensional problems.  Indeed,   almost always,  
one-dimensional  solvers for the Helmholtz operator  do  not translate well, and/or their results are not good predictors for two- or three dimensional solvers  because  the different level of difficulty in higher dimensions.  The main challenge there s is a richness and a high oscillatory character of the near-kernel of the Helmholtz operator, the problem that only very modestly presents itself in one dimension. 

Our position is very different: the  existing wave-ray algorithm successfully  deals (in two dimensions) with the near-kernel components. Its weakness  is  its limited applicability:  to this date, it has been developed only for constant wave numbers.   Therefore, our goal  is not to overcome the main challenge that is not fully present in 1D, but rather to extend the approach to variable, including discontinuous numbers,  dealing with  phenomenon  that   appear in all  dimensions.
\end{remark}

The  Helmholtz operator (\ref{eq:helmholtz}), especially  with constant $k(x) \equiv k$,  is expected   to be as   easily solvable  as   the  Laplace operator,  a poster child for  a successful   application of multigrid ideas. To the contrary,  Helmholtz equations is a completely different story. The main, though not the only one,  challenge in solving the Helmholtz equation iteratively, in particular using multigrid, is the set of its  near-kernel components ({\sl nkc}), the ones that satisfy   
 \be L v \approx 0
 \ee
 or, in a discrete formulation, 
\be
A^h v^h \approx 0^h.
\label{eq:nkc}
\ee
They play a significant role in any multigrid solver, including the ones described here.
In one dimension,  for instance, components (\ref{eq:nkc})  in the interior of $\Omega$  are of the form $ v(x)  = e^{i \omega x}$, $|\omega| \approx  k$. 
They  are too oscillatory to be accurately approximated on the coarse (the coarsest, for an adequate multigrid efficiency) scale,  and they have too small {\sl  relative residuals}  to be efficiently treated on fine grids (details can be found e.g., in \cite{BL97}).  Thus, these components is  a  liability for  most multigrid solvers as  they are virtually untreatable by standard multigrid techniques that involve using directly the Helmholtz operator,  e.g., \cite{elman}  or its shifted complex modification \cite{vuik1, vuik2}.  Such solvers  therefore can be used only  as  preconditioners to  Krylov  methods,  leaving to the latter  the  task of  eliminating (\ref{eq:error}).  This  strategy  leaves no chances for a scalability -- 
independence of the algorithm's performance  on the parameters of the problem and of the  solver. The  lack of scalability is especially pronounced in higher (than one) dimensions where the number of near kernel components increases  for larger $k$ and smaller $kh$:   such algorithms 
 can be only used for a limited range of wave numbers and using discretization with a mediocre resolution. \\
 
Another challenges include   divergence of {\sl linear} relaxation routines, such as the Gauss-Seidel, the Jacobi, and the SOR schemes, when applied to 
(\ref{eq:discrete}) and large phase discretization errors for components (\ref{eq:nkc}) .  
The latter can be improved, as often used  in  one-grid algorithms,  by employing  high-order discretization schemes. 

The wave-ray ({\sl wr}) approach, starting with  the original \cite{Thesis},  is based on  a different philosophy  on every one of these issues . First, at each multigrid iteration, it employs a special  treatment
for the {\sl nkc} that  efficiently reduces all of them, independently on their number and variety, serving as an actual solver rather than a  preconditioner.  Second, the algorithm considers   the original Helmholtz equation  and employs standard relaxation schemes.
Finally, since the algorithm shows a (nearly) optimal, aka linear, dependence on the problem's size, i.e.,  there is no need keep the size very small (for instance small enough for direct methods) and to employ high-order discretization schemes.  (The costs for using  such schemes in  the multigrid framework is not limited to the increased costs on the finest grid  -- to benefit from a high-order discretization, one must also use expensive high-order (larger stencils) interpolation and coarse-grid operators.  In case of Algebraic multigrid, for example,  this  leads to increasingly dense coarse grid operators, prohibiting using a fully efficient MG.)   

The  {\sl wr} algorithm  showed to be efficient for the most difficult, numerically,  Helmholtz operator  -- the one with  constant wave numbers, both in one \cite{Thesis} and two dimensions, \cite{BL97, BL06}.  Our goal is to extend its applicability to problems with   variable wave numbers, both continuous and discontinuous, while preserving all  special  features regarding the near-kernel treatment.  Here we propose strategies,  albeit implemented in one dimension, of how to achieve that, but with higher dimensions in mind -- all strategies are clearly  expendable to higher dimensions. 

The rest of the paper is organized as follows. In Section \ref{sec:gm}, a brief explanation of the wave-ray idea  is given,  followed by description of the wave-ray algorithm, its  parameters and numerical results for the original solver,  Section \ref{sec:parameters}. In Section \ref{sec:amg}, a new algebraic multigrid (AMG) version is introduced, accompanied with numerical results for constant wave numbers.
Modifications needed to apply the {\sl wr} approach, both geometric and algebraic,    for  continuous $k(x)$,  resulting in numerical experiments  presented  in Section \ref{sec:variable};  versions used for discontinuous wave numbers and the corresponding numerics, appear in Section \ref{sec:jump}.
 A  brief discussion of higher dimensions is given in Section \ref{sec:conclusions}.

\section{Geometric Multigrid (GMG)}
\label{sec:gm}
  The {\sl wr} approach is  based on  two observations.  First,  a standard  multigrid applied to the Helmholtz equation  efficiently  reduces all but the near-kernel error components.  Second,  the  dominant part of  unreduced  components can be represented in the form 
\begin{equation} e(x) = \han (x) \en + \hap(x) \ep,
\label{eq:exps}
\end{equation}
where  each function $\ha_\pm(x)$ is smooth, compared to the exponents. 
\begin{remark}
 It is important to note that (\ref{eq:exps}) is not  used to represent a solution, generated by some right-hand-side, but rather an  unreduced  error which mostly depends on $L$ (and $A^h$) as it largely consists of its  low energy modes. 
 Representation (\ref{eq:exps}) is rich:  It includes, for example,
 all exponents with frequencies $ \omega$,  $0 < |\omega| < 2k$, and many other functions, for instance the Hankel functions away from the origin.   As  analysis shows and numerical experiments confirm, \cite{Shifted},  in practice,  the actual range is much smaller: $( 1- \alpha_1) k \le |\omega| \le (1+\alpha_2)k$, with $\alpha_1 \approx \alpha_2 \approx 0.3$, depending on the relaxation regiment employed.
\end{remark}

In the {\it wr} approach, the task of computing the  oscillatory $e(x)$ is reduced to approximating two smooth   functions $\ha_\pm(x)$.  Given (\ref{eq:exps}),  the residual corresponding to $e$: $r = f -  A e$,  has a similar approximate representation:
 \begin{equation} r(x) = \hrn (x) \en + \hrp  \ep,
 \label{eq:res:exps}
 \end{equation}
 with smooth  residual functions $\hr_{\pm}$.  (Representation (\ref{eq:res:exps}) is  less accurate than (\ref{eq:exps}) as higher energy error (oscillatory) components with small amplitudes in $e(x)$ have larger amplitudes in $r(x)$.) 
  
Approximation  of  $\ha_\pm$ numerically requires two main ingredients: 
  \begin{itemize}
 \item Discrete coarse-grid,  operators;
  \item Residuals $\hr_\pm$ defined on the same coarse  grid. 
  \end{itemize}

  \begin{remark}
   As a nod to geometric optics terminology, all smooth hat functions are called {\sl } ray functions, the equations that describe them -- {\sl ray} equations, and the grids on which they are represented in the discrete formulation, {\sl ray} grids. Similarly, everything that is related to the Helmholtz part of the solver is called {\sl wave} functions, {\sl wave } operators  and {\sl wave} grids.  
  \end{remark}
\subsection{Ray operators and ray residuals}
\label{subsec:exps}

In geometric multigrid  each  discrete operator is a  discretization of an underlined   differential operator that describes the unknown function in continuum.  Ray differential equations  are discovered by applying  differential  (\ref{eq:helmholtz}) to (\ref{eq:exps}), yielding 
\begin{equation}
L \, e(x) =  \en \hLn \han (x) + \ep\hLp \hap (x) = r = \en \hrn (x)  + \ep \hrp (x),
\end{equation}
where 
\begin{equation}
\hLn\han= \han''  - 2 ik \han' \quad \mbox{and} \quad  \hLp\hap= \hap''  + 2 ik \hap' .
\label{eq:raydif} 
\end{equation}
Similarly,   Sommerfeld boundary conditions are translated in terms of $\ha_{\pm}$, resulting in, at $x = a$
\be
\en \han' = 0, \quad \ep [\hap'+2ik\hap] =0,
\ee
and, at $x=b$, 
\be
\en[\han'-2ik \han] =0, \quad \ep \hap' =0. 
\ee

Advantage of using exponential basic functions $\exp(\pm ikx)$ here is clear:  the high oscillations in (\ref{eq:exps}) are completely removed from  the ray description, and the two ray equations can be separated into two individual systems:  
\begin{equation}
\left\{ \begin{array}{lll}
         \hLn \han = \hrn , & \han'(a) = 0,  & [\han'-2ik \han](b) = 0;\\
         \hLp \hap = \hrp,  &[\hap'+2ik\hap](a) =0, & \hap'(b) =0.
         \end{array}  \right. 
\end{equation}
The ray operators are  discretized on the {\sl ray} grid with the mesh-size that satisfies $kH \approx \pi$, using a four-point stencil defined on staggered grids (details and motivation can be found in \cite{Thesis, BL97} ).
Ray  residuals are also approximated on  scale $H$, using a  {\sl separation} procedure  that  relies on the following properties.  For a smooth  function $g^h$ and exponential functions  $e^{\pm 2ikx }$, both  defined on   fine scale, $h$, holds
\be || g^h - P_H^h R_h^H g^h || \ll || g^h|| \quad \mbox{and}  \quad R_h^H (g^h e^{\mp 2ikx }) \approx {{ 0}^H}, 
\label{eq:orth}
\ee
where $P^h_H$ and  $R_h^H$ are a linear interpolation  and  a full weighting acting from  $H$ to $h$ and from $h$ to   $H$,  respectively.  
 In other words, $R^h_H$ preserves smooth functions and nearly eliminates functions close, in the frequency space,  to $ e^{\pm 2ikx }$. These  considerations lead  approximation procedure for $\hr_{\pm}^H$. Given  the finest-grid wave residual, $r^h$, ray residuals can be computed as  
\[  R_h^H (\ep r^h)   \approx R_h^H(\hrn^h + e^{2ikx} \hrp^h)  \approx \hrn^H. \]
and 
\[  R_h^H (\en r^h)   \approx R_h^H(e^{-2ikx}\hrn^h + \hrp^h)  \approx \hrp^H. \]
\section{WR algorithm}
\label{sec:parameters}
The wave-ray algorithm consists of two parts: a standard V-cycle applied to the Helmholtz equation and the additional correction by the two coarse-grid ray systems,  as  described in the following  Pseudocode. \\

\noindent 
\begin{tabbing} 
\quad  \= \qquad  \= \qquad \= \qquad \= \kill

\underline{ Wave Cycle: $ WaveCycle(\ell, r^\ell)$}\\ \\
The input: Current grid  $\ell$,   the residual $r^\ell$\\
The output: New finest grid correction  $e^1$ for $\ell = 1$ \\ 
The MG framework: {\sl Correction Scheme} \\
\>if $\ell == L$ \\
\>$\%$ On the coarsest grid $L$ \\
\>\>$e^\ell = Relaxation(A^\ell, r^\ell)$ (or $e^\ell = (A^\ell)^{-1} r^\ell$)\quad $A^\ell$ is a finite difference approximation  (\ref{eq:discrete}) on scale $\ell$\\
\>\>  $e^{\ell-1} = r^{\ell-1} + P_{\ell}^{\ell-1} e^\ell$ \quad  ($  P_{\ell}^{\ell-1} $ is a linear interpolation from grid $\ell$ to grid $\ell-1$)\\
\>else\\
\>\>$\%$ On grids $\ell = 1, \dots, L-1$ \\
\>\>$e^\ell = Relaxation(r^\ell)$\\
\>\> $ r^{\ell+1} = R_{\ell}^{\ell+1} (r^\ell-A^\ell e^\ell)$  \quad  ($  R^{\ell+1}_{\ell} $ is a full weighting from grid $\ell$ to grid $\ell+1$)\\
\>\> $e^\ell = e^\ell + WaveCycle(\ell+1, e^{\ell+1})$ \\
\>\>$e^\ell = e^{\ell} + Relaxation(r^\ell-A^\ell e^\ell )$\\
\>\> if $\ell > 1$\\
\>\>\> $e^{\ell-1} = e^{\ell-1} + P_{\ell}^{\ell-1} e^\ell$ \\
\>\>end if \\
\>end if\\ \\

\underline{ Ray Cycle: $ RayCycle(r)$}\\ \\
The input:  Finest grid wave residual $r \equiv  r^1$\\
The output:  New finest grid correction   $e \equiv e^1$\\
The MG framework: {\sl Correction Scheme} \\
\>$ [\hrn^H, \hrp^H ] = Separation(r^) $ \\
\> Solve $\hLn \han = \hrn$ and $\hLp \hap = \hrp$ (by Gauss-Seidel relaxation in the propagation (negative and positive) direction);\\
\> $e  = \en (P_H^h \han^H) + \ep (P^h_H \hap^H)$\\ \\

\underline{ Wave Ray Cycle: $ WaveRayCycle(A, b, \hat{x})$}\\ \\
The input:  $A \equiv A^1$, $b \equiv b^1$, the  current approximation $\hat{x} \equiv x^1$\\
The output:  New finest grid approximation $\tilde{\tilde{x}}$\\
The MG framework: {\sl Correction Scheme} \\
\> $\tilde{x} = \hat{x}+WaveCycle(b-A \hat{x})$\\
\> $\tilde{\tilde{x}} = \tilde{x} + RayCycle(b-A \tilde{x})$\\
\quad  \= \qquad  \= \qquad \= \qquad \= \kill
\end{tabbing}

All numerical methods in this paper are run for the wave-ray cycle (and its variation described in the upcoming Sections) with the fixed set of parameters 
that depend on values  
\[  k_{max} = \max_{x\in \Omega}  \, k(x), \quad k_{min} = \min_{x\in \Omega} \, k(x).
\]
They are 
\begin{itemize}
\item The finest mesh-size satisfies $k_{max} h <  2\pi/10$, typically $k_{max} h \approx 0.3$ or less;
\item The ray mesh-size satisfies $ \pi/2 <  k_{min} H \le \pi$; 
\item In the wave cycle the Gauss-Seidel relaxation is used on all grids except the ones with the intermediate  mesh-sizes $\tilde{h}$: $\pi/4 < k_{max} \tilde{h} \le \pi/2$;
\item The number of relaxation in the wave cycle: one  pre- and  post relaxation if Gauss-Seidel is employed; two pre- and post sweeps for Kaczmarz relaxation; 
\item The ray relaxation is the Gauss-Seidel that in this case is an almost direct solver; 
\item The number of relaxation in the ray cycle: two  per ray component. 
\end{itemize}

\noindent
In all Tables,  computations are performed until the residual satisfies: 
\be \frac{||r_m||}{||r_0||} < 10^{-6}, 
\label{eq:resid_reduction}
\ee
where $||r_0||$ and $||r_m||$ are the $\ell_2$ norms of the initial residual and the residual after $m$ wave-ray cycles.  
The numerical results for the Helmholtz equation (\ref{eq:helmholtz}) with constant $k$ are presented in Table \ref{tab:gm:const}.

\begin{table}[htp]
\centering{
\begin{tabular}{l | c  c c r}
\hline 
$k$ & 40  & 80 & 160 & 320\\ \hline
$kh = 0.625$  & 13 & 14  & 27  &  14 \\ \hline
$kh = 0.3125$  & 12 & 12  & 13  &  14 \\ \hline
$kh = 0.15625$  & 15 & 15  & 16  &  17 \\ \hline
\end{tabular}}
\vspace{0.5cm}
\caption{The number of gmg-WR  cycles needed to satisfy (\ref{eq:resid_reduction}).  
The results are computed for different choices of $k$ and the finest  $h$; the number of levels varies from $L = 5$ ($k = 40, kh = 0.625$) to $L = 10$ ($k=320, kh= 0.15625$).}
\label{tab:gm:const}
\end{table}
The results will serve as a  benchmark for  all other variants of the algorithm, developed and applied to problems with different types of wave numbers.  
  
  \begin{remark}
  Certain values of $k$   have to be chosen for  numerical experiments; These values used  are {\sl  not } the best case scenarios   for which the WR algorithm performs the best. The results slightly depend  values of $k$ and $kh$. The difference is due to different smoothing rates of relaxation schemes  and, more importantly, to accuracy of the separation procedure. The algorithm can be tailored  to accommodate the separation to  a particular value of $kH$;  the   full weighting employed here works most accurately for  $kH = \pi$ - and such values do no appear in the paper. 
  \end{remark}
\section{Algebraic multigrid wave-ray (amgWR) algorithm }
\label{sec:amg}

The geometric wave-ray algorithm relies on the knowledge of  analytical near-kernel components of the differential Helmholtz operator $L$. An  alternative {\sl algebraic} version  uses numerical approximations instead of analytical functions and proceeds to directly compute discrete coarse grid operators, both in wave and ray representation, bypassing the differential ones.  The advantage of the approach is its ability to numerically adjust the basic functions to better satisfy (\ref{eq:nkc}) .

The amgWR  algorithm consists of two parts: wave and ray. In both,  Galerkin approach is used to form a coarse grid operator from its finer grid predecessor
\be
A^{h_c} = (P_{h_c}^{h_f})^{t} A^{h_f} P_{h_c}^{h_f}.
\label{eq:galerkin}
\ee 
Here  $A^{h_f}$ is a given fine-grid operator,  $P_{h_c}^{h_f}$ and $(P_{h_c}^{h_f})^{t}$ are  a linear interpolation  and its transposed,   scale  $h_c$ to a fine scale $h_f$, and {\sl t } stands for {\sl transposed}, $h_c > h_f$;  $A^{h_c}$ is the resulting coarse grid operator.

\subsection{Wave part}
\label{subsec:wave}

The main decision  in the AMG  is a  choice  of prolongation operators,$P$. Typically, they are constructed  to accurately transfer  the near-kernel components of the finest-grid operator: in the Helmholtz case these are components  (\ref{eq:exps}).  In the WR approach,  however, the wave part is not responsible for their treatment, and, therefore, $P$ is not defined by their character.  Two remaining types of components: highly oscillatory and physically smooth  are well served by a standard polynomial interpolation, similarly to  processing of Laplace operator.  Indeed, the former are treated on sufficiently fine grid where they remain relatively smooth, while the later are reduced to the coarsest grids from which they are accurately interpolated due to their smoothness.  All coarse-grid  wave operators, $\ell =2, \dots, L$ 
are computed using (\ref{eq:galerkin}) , with $A^1 \equiv A^h$.  
\begin{remark}
Operators $A^\ell, \quad \ell =2, \dots, L$ are very similar to the one of the discrete geometric operators -- they too are some discretizations of $L$  and can be easily substituted for one another.
\end{remark}
\subsection{Ray part}
\label{subsec:ray}

Construction of  coarse grid ray operators starts on the finest grid, using the  Helmholtz discrete operator, $A^h$.  Using  $\en$ and $\ep$, substitution of 
 error  in the form (\ref{eq:exps}  into the wave residual equation $A^h e^h = r^h$, produces finest-grid ray operators  $\hat{A}^h_{-}$ and $\hat{A}^h_+$:
\be A^h e^h = \en \hat{A}^h_{-} \han^h + \ep \hat{A}^h_{+} \hap^h, \ee
with their  stencils   given by
\[   \hat{S}^h_{-} = \frac{1}{h^2} \biggl [ e^{ikh} \quad (-2 +k^2h^2) \quad e^{-ikh} \biggr ] \]
and 
\[  \hat{S}^h_{+} =  \frac{1}{h^2} \biggl [e^{-ikh}\quad  (-2+k^2h^2) \quad  e^{ikh} \biggr]. \]
Ray operator $\hat{A}^H_{\pm}$ on scale  $H$ are then computed as 
\be \hat{A}_{\pm}^{H} = (P^h_{H})^t (\hat{A}_{\pm}^h) (P^h_{H}). 
\label{eq:rayH}
\ee 
A   linear interpolation again is sufficiently accurate for smooth $\ha_\pm$. 

\noindent
Ray residuals are approximated by 
\[ \hrn^{H} = (P_{H}^h)^t (\ep r^h) \]
and 
\[ \hrp^{H} = (P_{H}^h)^t (\en r^h). \]
The {\it amgWR} algorithm is applied to   (\ref{eq:discrete}) with  constant $k$, Table \ref{tab:amg:const},  with parameters as in 	 Table \ref{tab:gm:const}.  Again, the results  clearly improve or at least stabilize when the algorithm is applied for small $kh$. This  effect becomes  more pronounced for larger wave numbers, they are similar to the ones of the {\it gmgWR}.

\begin{table}[htp]
\centering{
\begin{tabular}{l | c  c c r}
\hline 
$k$ & 40  & 80 & 160 & 320\\ \hline
$kh = 0.625$  &  16 & 34  & $>50$  & $> 50$ \\ \hline
$kh = 0.3125$  & 11  & 13  & 18   &   43  \\ \hline
$kh = 0.15625$  & 11 &  10  & 12  & 18 \\ \hline
$kh =  0.078125$  & 11 & 10  & 12  & 14 \\ \hline
\end{tabular}}

\vspace{0.5cm}
\caption{Presented here is the number of the AMG wave-ray cycles needed to satisfy (\ref{eq:resid_reduction}).  The results are obtained both for different values of $k$ and $kh$ where $h$ is the finest scale mesh-size.}
\label{tab:amg:const}
\end{table}

\section{Continuous wave numbers: {\it amgWR} and {\it amgWR(c)}}
\label{sec:variable}
As an example of a variable wave number here 
\be k(x) = k_0\sqrt{1 + m(x)}, \quad \mbox{with} \quad |m(x)| <  1, \ee
is considered, and in the numerics 
\be m(x) = \alpha \cos(\beta x)
\ee
is used. Such form allows to investigate how  amplitude and oscillations in $m(x)$ affect algorithmic performance.  First,  the {\it amgWR} with exponential basic functions  $u_\pm(x) = e^{\pm i k_0 x}$ is tested, Tables  \ref{tab:amg:var1}-  \ref{tab:amg:var2}. 
Table \ref{tab:amg:var1}  presents the results for smooth $k(x)$ with $\beta$ fixed at  $1$,  computed  for two different values of $k_0$. 

\begin{table}[htp]
\centering{
\begin{tabular}{l | c  c c r}
\hline 
$\alpha  $ & 0.1  & 0.2 & 0.4 & 0.8\\ \hline
$k_0 = 40$  & 10  &   10  &   11 &  13 \\ \hline
$k_0 = 160$  &  11 & 11   & 12  &  14 \\ \hline

\end{tabular}}

\vspace{0.5cm}
\caption{The number of the AMG WR cycles needed to satisfy (\ref{eq:resid_reduction}); $\beta = 1.$}
\label{tab:amg:var1}
\end{table}

In Table \ref{tab:amg:var2},    wave numbers are considered for a variety of $k_0$ and $\beta$ values and  $\alpha= 0.5$.  The results suggest that   $\alpha $ value  does not impact convergence,  at least when it remains bounded by one.  Increase  in  $\beta$ leads to slowdown and   divergence. 

\begin{table}[htp]
\centering{
\begin{tabular}{l | c  c c c  r}
\hline 
$\beta  $ & $0.1k_0$  & $0.25k0$ & $0.5k_0$ & $0.75 k_0$& $ k_0$   \\ \hline
$k_0 = 25$  &  12&   11  &  11 &  30   & 41  \\ \hline
$k_0 = 50$  & 12 & 11 &  11  &  D  & D \\ \hline
$ k_0 = 100$  & 13  & 13     &  13   &  D   &  D \\ \hline
$ k_0 = 200$  &  15 &  21  &  19  &  D  & D  \\ \hline
\end{tabular}}
\vspace{0.5cm}
\caption{The number of the AMG WR cycles needed to satisfy (\ref{eq:resid_reduction}); $beta$ is proportional to $k_0$, with $k_0/\beta$ varying from $0.1$  to $1$, $\alpha = 0.5$;  "D" stands for divergence.}
\label{tab:amg:var2}
\end{table}

The next  variation {\it amgWR(c)}   improves  convergence   by modifying the basis functions, while still maintaining directions as defined by the exponents.  The new basis functions are  sought in the  ray form 
\be u_{\pm}(x) = \hat{u}_{\pm}(x) e^{\pm ik_0 x}, 
\label{eq:ex:modified}
\ee
where  $\hat{u}_{\pm}$  satisfy  
\be  \Delta \hat{u}_{\pm} (x) \pm 2ik_0 \hat{u}_\pm + k_0^2 m(x) \hat{u}_{\pm}(x) = 0, 
\label{eq:uhat}
\ee
with boundary conditions derived from (\ref{eq:bc}), 
\begin{equation}
\left\{ 
\begin{array}{lr}
         \hat{u}'_\pm(a) + i(k(a) \pm  k_0) \hat{u}_\pm(a) = 0, & \\
           \hat{u}'_\pm(b) - i(k(b) \mp k_0) \hat{u}_\pm(b) = 0. &  \end{array} \right. 
\end{equation}

Discrete equations for $\hat{u}_{\pm}$  on scale $H$ can be obtained either geometrically or algebraically, and functions $\hat{u}_{\pm}$  are computed, as  discussed in Section{\ref{subsec:exps}  using $\hat{r}_{\pm} = 0$.  
 After $\hat{u}_{\pm}^H$ are approximated and interpolated to the finest scale, $h$, 
the modified  basis functions (\ref{eq:ex:modified}) are reconstructed and employed in the  wave-ray algorithm.  The  numerical results are given in Table \ref{tab:amg:var3}, using parameters identical to the ones in Table \ref{tab:amg:var2}.

\begin{table}[htp]
\centering{
\begin{tabular}{l| c  c c c  r}
\hline 
$\beta  $ & $0.1k_0$  & $0.25k0$ & $0.5k_0$ & $0.75 k_0$& $k_0$   \\ \hline
$k_0 = 25$  &  D &   D  &  13 &  13   & 18  \\ \hline
$k_0 = 50$  & D & D &  13  &  14  & 16 \\ \hline
$ k_0 = 100$  & D  & D     &  13   &  14  &  21 \\ \hline
$ k_0 = 200$  &  D &  D  &  17  &  14  & 20 \\ \hline
\end{tabular}}
\vspace{0.5cm}
\caption{The number of the {\it amgWR(c)} cycles needed to satisfy (\ref{eq:resid_reduction}); $beta$ is proportional to $k_0$, with $k_0/\beta$ varying from $0.1$
 to  $1$, $\alpha = 0.5$;  "D" stands for divergence.}
\label{tab:amg:var3}
\end{table}
The {\it amgWR(c)} modified approach compliments {\it amgWR} : it diverges when the former converges well and shows good results when the other diverges.  The combined Table \ref{tab:combined} gives  the best of the two approaches' results. Clearly, these are just the first steps in understanding what are the best strategies in dealing with oscillatory wave numbers -- but the results  are quite  encouraging.

\begin{table}[htp]
\centering{
\begin{tabular}{l| c  c c c  r}
\hline 
$\beta  $ & $0.1k_0$  & $0.25k0$ & $0.5k_0$ & $0.75 k_0$& $k_0$   \\ \hline
$k_0 = 25$  &  12 &   11  &  11 &  13   & 18  \\ \hline
$k_0 = 50$  & 12 & 11  &  11  &  14  & 16 \\ \hline
$ k_0 = 100$  & 13  & 13     &  13   &  14  &  21 \\ \hline
$ k_0 = 200$  &  15 &  21  &  17  &  14  & 20 \\ \hline
\end{tabular}}
\vspace{0.5cm}
\caption{The number of the {\it amgWR} or {\it amgWR(c)} cycles, whichever is the smallest,  needed to satisfy (\ref{eq:resid_reduction}); $beta$ is proportional to $k_0$, with $k_0/\beta$ varying from $0.1$
 to  $1$, $\alpha = 0.5$.}
\label{tab:combined}
\end{table}

 
\section{Discontinuous Wave Numbers: {\it amgWR(d)} and {\it gmgWR(d)}}
\label{sec:jump}
Considered here is a discontinuous  wave number in  the form
\begin{equation}
k(x)  = \left\{ \begin{array}{ll}
         k_1, & \mbox{if $x \leq \bx $},\\
        k_2, & \mbox{if $x >  \bx $},\end{array} \right. 
        \label{eq:jumpk}
\end{equation}
with   $k_1 \ge k_2$ for definiteness.
First, the {\it amgWR}    is  applied using basic functions
\begin{equation}
u_\pm (x)  = \left\{ \begin{array}{ll}
         e^{\pm i k_1 x}, & \mbox{if $x \leq \bx $},\\
        e^{\pm i k_2 x }, & \mbox{if $x >  \bx $}.\end{array} \right. 
        \label{eq:jumpk:amg}
\end{equation}
The results for $k_1 = k$ and $k_2 = \gamma k$, $0 < \gamma < 1$  are given in Table \ref{tab:disc:amg}, they  are restricted to cases when  both $k_1$ and $k_2$ are  large enough 
  to benefit from  ray representation.  The algorithm can be easily adjusted to accommodate small values of $k$ but it is not in the scope of  this paper. \\

\begin{table}[htb]
\centering{
\begin{tabular}{l | c  c c r}
\hline 
$k$ & 40  & 80 & 160 & 320\\ \hline
$\gamma = 0.8$  & 18 (13) & 31(14)  & 36(19) & $>50$(20)  \\ \hline
$\gamma = 0.5$  & 18(13)  & 33(15) & 37(21)& $>50$(23)\\ \hline
$\gamma = 0.25$  & 18(13) & 33(15) & 37(21) & $>50$ (23)\\ \hline
\end{tabular}}
\vspace{0.5cm}
\caption{The number of the {\it amgWR} cycles needed to satisfy  (\ref{eq:resid_reduction}):  with  basis functions (\ref{eq:vp})-(\ref{eq:vm}) and (in parenthesis) with  each of  (\ref{eq:vm})-(\ref{eq:vp}) functions  pre-smoothed by one wave cycle.}
\label{tab:disc:amg}
\end{table}
There is  a clear benefit of   {\sl pre-smoothing} of the basis functions  which are discontinuous  at $\bx$;  it is very likely that  a local processing  near $x = \bx$ would be sufficient, optimization of the  pre-smoothing strategy is a subject of future investigation.  Overall,  the algorithm loses  efficiency compared to constant $k$.

\subsection{Residual Separation for Discontinuous Wave Numbers}
\label{subsec:separation}
Ray residuals $\hr_{\pm}$ are approximated on scale that satisfies $\pi/2 \le kH \le \pi$.  Wave numbers (\ref{eq:jumpk})  yield two different scales:
$H_1$, with $\pi/2 \le k_1 H_1 \le \pi$, and $H_2$, with $\pi/2 \le k_2 H_2 \le \pi$. If $\gamma > 1/2$  then $H_1 = H_2 = H$, and the standard separation procedure described in Section \ref{sec:amg} works perfectly well.  However,  if $\gamma \le 1/2$ then $H_2 \ge  2 H_1$ and the separation and possibly the approximation of $\hanp$ is better done on different scales   for $x \le \bx$ and $x > \bx$.

Current  implementation is a compromise: ray functions are approximated on the same scale $H = \min\{H_1, H_2\} = H_1$  throughout $\Omega$, meaning that both residuals $\hr_\pm$ are eventually assembled there. However,  scale $H_2$ is used in the separation routine in the following way.
\begin{itemize}
\item Residuals $\tilde{r}^H_\pm$ are  approximated as discussed in Section \ref{sec:amg}; 
\item For $x > \bx$ the averaging continuous to scale $H_2$: 
\[  \biggl(\tilde{r}^{H_2}_\pm\biggr)_{x > \bx} =  \biggl( R_{H}^{H_2} \tilde{r}^H_\pm\biggr)_{x > \bx} 
\]
Thus, the residuals are properly separated on $x > \bx$; 
\item Ray residuals are reconstructed on the ray grid
\begin{equation}
\hr_\pm^H  = \left\{ \begin{array}{ll}
         \tilde{r}^H_\pm, & \mbox{if $x \leq \bx $},\\
        P_{H_2}^H \tilde{r}^{H_2}_\pm, & \mbox{if $x >  \bx $},\end{array} \right. 
        \label{eq:sepresid}
\end{equation}
where $P_{H_2}^H$ is a linear interpolation from scale $H_2$ to scale $H$; it works because 
$\tilde{r}^{H_2}_\pm$ are smooth.
\end{itemize}

\subsection{AMG and  Geometric Optics}

The next modification  is based  the  laws of geometric optics. 
If an  incident   wave is given by $e^{-ik_2 x}$, propagating in the negative direction and entering $\Omega$ from the right, solution $L\overline{u}_{-} = 0$ is of  the form
\begin{equation}
\overline{u}_{-}(x)  = \left\{ \begin{array}{ll}
     C^t_{-} e^{-ik_1 x} ,    & \mbox{if $x \leq \bx $},\\
    e^{-ik_2 x} + C^r_{-} e^{ik_2 x},    & \mbox{if $x >  \bx$},\end{array} \right.
        \label{eq:itrn}
\end{equation} 
with  reflection and transmission coefficients $C_{-}^r$ and $C_{-}^t$ given by:
\begin{equation}
C^r_{-} = \frac{k_1-k_2}{k_1+k_2}\, e^{-2ik_2\bx}, \quad C^t_{-}= \frac{2k_2}{k_1 + k_2}\, e^{i(k_1-k_2)\bx}.
\end{equation}
For $\overline{u}_{+}(x)$, with an incident  wave $e^{ik_1 x}$, entering $\Omega$ from  the left, the solution is given by 
\begin{equation}
\overline{u}_{+}(x)  = \left\{ \begin{array}{ll}
     e^{ik_1 x} + C^r_{+} e^{-ik_1 x}, & \mbox{if $x \leq \bx $},\\
     C^t_{+} e^{ik_2 x}, & \mbox{if $x >  \bx$}.\end{array} \right.
        \label{eq:itr}
\end{equation} 
with   coefficients $C_{+}^r$ and $C_{+}^t$ given by:
\begin{equation}
C_{+}^r = \frac{k_1-k_2}{k_1+k_2}\, e^{2ik_1\bx}, \quad C_{+}^t = \frac{2k_1}{k_1 + k_2}\, e^{i(k_1-k_2)\bx}.
\end{equation}

\noindent
The new basic functions are composed of  components  of (\ref{eq:itrn})-(\ref{eq:itr})  that propagate in the same direction: 
  \begin{equation}
 u_{-}(x)  = \left\{ \begin{array}{ll}
       C^t_{-}  e^{-ik_1 x}, & \mbox{if $x \leq \bx $},\\
         e^{-ik_2 x}, & \mbox{if $x >  \bx $}.\end{array} \right.
          \label{eq:vm}
  \end{equation}
  and
 \begin{equation}
 u_{+}(x)  = \left\{ \begin{array}{ll}
       e^{ik_1 x}, & \mbox{if $x \leq \bx $},\\
        C^t_{+} e^{ik_2 x}, & \mbox{if $x >   \bx$}.\end{array} \right.
         \label{eq:vp}
 \end{equation} 
Table \ref{tab:disc:amg2}  presents numerical results for different values and ratios of $k_1$ and $k_2$.  Once more,  the results significantly improve when the basic functions  are preprocessed and the efficiency becomes  similar to the one obtained for constant wave numbers. 
\begin{table}[htb]
\centering{
\begin{tabular}{l | c  c c r}
\hline 
$k $ & 40  & 80 & 160 & 320\\ \hline
$\gamma = 0.8,( p = 0)$  & 30 (11) & 27(21)  & 18(12) & 13(12)  \\ \hline
$\gamma = 0.5, (p = 1)$  & 45(11)  & 44(11) & 31 (13)& 36 (15)\\ \hline
$\gamma = 0.25,(p = 2)$  & 23(12) & 23(12) &  33(13) &  27(14)\\ \hline
\end{tabular}}
\vspace{0.5cm}
\caption{ The number of the {\it amgWR(d)}  cycles needed to satisfy (\ref{eq:resid_reduction}), with the basis functions defined by (\ref{eq:vm})-(\ref{eq:vp}) without and with (in parenthesis)  pre-smoothing by one  AMG wave cycle;  parameter $p$ describes the ratio between the ray scales $H_2 = 2^p H_1$ in Section \ref{subsec:separation}.}
\label{tab:disc:amg2}
\end{table}

\subsection{{\it gmgWR(d)}:  Adaptation of  Ray Operators}

The original {\it gmgWR} has a limited application for Helmholtz operators with discontinuous wave numbers. As shown in  Table  \ref{tab:disc:amg2}, its   performance deteriorates even for $\gamma \approx 1$, and it completely falls apart as $\gamma$ becomes smaller.  
\begin{table}[htb]
\centering{
\begin{tabular}{l | c  c c r}
\hline 
$k$ & 40  & 80 & 160 & 320\\ \hline
$\gamma = 0.95$  &  21 &  55  & 51  & 43\\ \hline
$\gamma = 0.8$  &  52 &  45  & 80  & 68\\ \hline
$\gamma = 0.4$  &  24  & 62  & D  & D \\ \hline
$\gamma = 0.25$  & D  & D  &  D &  D \\ \hline
\end{tabular}}
\vspace{0.5cm}

\caption{The number of the GMG WR cycles needed to satisfy (\ref{eq:resid_reduction}); D stands for {\sl divergence}. }
\label{tab:disc:gmg}
\end{table}

The reason for such a poor performance is a non-adequate description of  interfaces between  different media  by ray equations.  
To address that,  ray equation(s), with stencils that are crossed by such interfaces are modified  based on  information from a much more detailed finest wave grid. If, for example, 
$k(x-h) = k(x) = k_1$ and $k(x+h) =k_2$ and assuming 
    \begin{equation}
    u(x)  = \left\{ \begin{array}{ll}
             \ha_{\pm} e^{\pm ik_1 x}& \mbox{if $x \leq \bar{x}$},\\
             \ha_{\pm} e^{\pm ik_2 x} & \mbox{if $x >  \bar{x}$}, \end{array} \right. 
    \end{equation}
    then  substituting   $u(x)$ into  a discrete Helmholtz equation, centered at $\bx$,  and  subsequently   applying  Taylor expansion,  yields   a modified differential ray operator centered at $\bx$: 
\be
\begin{split}
& \frac{\hanp(x-h) e^{\pm ik_1(\bar{x}-h)} -2\hanp(\bar{x}) e^{\pm ik_1 \bar{x}}+\hanp(\bar{x}+h) e^{\pm ik_2(\bar{x}+h)}}{h^2} + k_1^2 \hanp(\bar{x}) e^{\pm ik_1 \bar{x}} = \\
& \hat{L}_{\pm} \hanp(\bx)+  \biggl[\hanp(x) + h\hanp'(x) + \frac{h^2}{2}\hanp''(\bx)\biggr]\frac{V_{\pm}}{h^2},
\end{split}
\label{eq:gmcor}
   \ee
   
where \[ V_{\pm} = {e^{\pm i(k_2-k_1)\bx}e^{\pm ik_2 h}-e^{\pm k_1 h}}.\]
Table \ref{tab:disc:gmg1} shows that modification of  one ray discrete equation for each of $\hanp$ leads to a significant improvement.

\begin{table}[htb]
\centering{
\begin{tabular}{l | c  c c r}
\hline 
$k$ & 40  & 80 & 160 & 320\\ \hline
$\gamma = 0.8$  &  18 &  21 & 23 & 23 \\ \hline
$\gamma = 0.4$  &  22 &  20 & 21 & 21 \\ \hline
$\gamma = 0.25$  & 21  & 25 & 25 & 21 \\ \hline
\end{tabular}}
\vspace{0.5cm}

\caption{The number of the {\it gmgWR(d)}  cycles needed to satisfy (\ref{eq:resid_reduction})}
\label{tab:disc:gmg1}
\end{table}
Although the results by the {\it gmgWR(d)}  are  modes compared to the  ones of the  {\it amgWR(d)}  , its  setup costs are smaller, and this could play a  role in higher dimensions, where   fast convergence  might be outweighed by  high setup costs.
Which approach will have a better overall performance in two- and three dimensions remains to be seen.

\section{Conclusions and extension to higher dimensions}
\label{sec:conclusions}
In this paper, steps  to extend the existing  geometric wave-ray algorithm  for the Helmholtz operator with  constant wave numbers to problems with both continuous and discontinuous wave numbers, are outlined.  The results, achieved  by modification of the  existing geometric and new  algebraic  versions of the wave ray algorithms, are comparable to the results shown by the original wave ray algorithm for constant $k$.  The summary of the results observed for different versions of the algorithm for various  types of the wave numbers are shown in Table \ref{tab:all}.
\begin{table}[htb]
\centering{
\begin{tabular}{l | c  c c r}
\hline 
$k$ & 40  & 80 & 160 & 320\\ \hline
Constant $k$: GMG  &  12 & 12  &  13  & 14    \\ \hline
Constant $k$: AMG   & 11(11)  & 13(10)  & 18 (12)  &  43 (14)\\ \hline
Continuous $k$ ($\alpha = 0.4, \beta = 20$): AMG  & 13 & 7 & 8   & 9\\ \hline
Discontinuous $k$ ($\gamma = 0.25$) GMG(d)  &  21 &  25 & 25 & 21 \\ \hline
Discontinuous $k$ ($\gamma = 0.25$) AMG(d)  & 14  & 14 & 15 & 16 \\ \hline
\end{tabular}}
\vspace{0.5cm}

\caption{The number of different types of  WR  cycles, generic   or adapted  to a particular type of wave numbers,  needed to satisfy (\ref{eq:resid_reduction}). For constant $k$ for the AMG WR the number in parenthesis is the number of cycles if the problem is considered on finer grid with $kh = 0.078125$. }
\label{tab:all}
\end{table}

The  extension to higher dimensions relies on an  analogue of the error representation (\ref{eq:exps}):
\[ e({\bf x}) = \sum_{\kappa =1}^K \hat{a}_\kappa ({\bf x}) e^{i{\bf k}^\kappa {\bf x}}, \quad {\bf x} \in \Omega \subset {\bf R}^d, \quad {\bf k}_\kappa \in {\bf R}^d, \quad |{\bf k}_\kappa| = k, \quad d  = 2, 3,\]
along with the  used in the original two-dimensional  {\it gmgWR(d)}  algorithm.  The frequencies ${\bf k}_\kappa$ are uniformly distributed along a circle (sphere) of radius $k$. The number of the basis functions $K = O(1)$ grows for higher dimensions;  for instance, for the two-dimensional solver in \cite{BL97} $K = 8$ was sufficient.

The  next step  is to extend  the ideas  presented in Sections \ref{sec:amg}-\ref{sec:jump}  to two dimensions.  The first results    in this direction were obtained for the two-dimensional {\it amgWR} solver, and they are promising \cite{NewWork}.  The new algorithm along with the original two-dimensional 
{\it gmgWR} \cite{BL97}  will  serve as a foundation for implementing  other strategies presented in this paper.  
\bibliographystyle{siam}
\bibliography{helm2}	
\end{document}